\documentclass[a4paper,10pt]{article}
\usepackage{amsmath}
\usepackage{amssymb}
\usepackage{bbm}
\usepackage[pdftex]{graphicx}

\title{The optimal constants in Khintchine's inequality for the case $2<p<3$\footnote{Supported by DFG project KO 962/10-1}}
\author{Olaf Mordhorst}

\begin{document}

\maketitle
\newtheorem{definition1}{Definition}[section]
\newtheorem{remark1}[definition1]{Remark}
\newtheorem{lemma1}[definition1]{Lemma}
\newtheorem{proposition1}[definition1]{Proposition}
\newtheorem{theorem1}[definition1]{Theorem}
\newtheorem{corollary1}[definition1]{Corollary}
\newtheorem{example1}[definition1]{Example}

\begin{abstract}
A mean step in Haagerup's proof for the optimal constants in Khintchine's inequality is to show integral inequalities 
of type \(\int(g^s-f^s)\mathrm{d}\mu\geq 0\). In \cite{NazarovPodkorytov} F.L. Nazarov and A.N. Podkorytov made Haagerup's proof much more clearer 
for the case \(0<p<2\) by using a lemma on distribution functions. In this article we want to treat the case \(2<p<3\)
with their technique.  
\end{abstract}

\section{Introduction}
We will always denote by \(\varepsilon_1, \varepsilon_2,\varepsilon_3, \dots\) a sequence of i.i.d. Bernoulli variables
on a probability space \((\Omega, \mathrm{P})\) with 
\(\mathbbm{P}(\{\varepsilon_1=-1\})=\mathbbm{P}(\{\varepsilon_1=1\})\) and we denote by \(\mathbbm{E}\) the expectation.
The Khintchine inequality states the following:
\begin{theorem1}\label{MainTheorem}
For every \(p\in\mathbbm{R}_{>0}\) there exist constants \(A, B >0\) such that for every \(n\in\mathbbm{N}\)
and every sequence \((a_k)\in \mathbbm{R}^n\) we have
\begin{align}
 A\|(a_k)\|_2\leq \left(\mathbbm{E}|\sum_{k=1}^n a_k\varepsilon_k|^p\right)^{1/p}\leq B \|(a_k)\|_2\label{KInequality}
\end{align}
The optimal constants for which (\ref{KInequality}) holds are 
\(A_p=\min\{1,2^{\frac{1}{2}-\frac{1}{p}}, 2^{\frac{1}{2}}\left(\frac{\Gamma(\frac{p+1}{2})}{\sqrt{\pi}}\right)^{\frac{1}{p}}\}\)
and \(B_p=\max\{1, 2^{\frac{1}{2}}\left(\frac{\Gamma(\frac{p+1}{2})}{\sqrt{\pi}}\right)^{\frac{1}{p}}\}\).
\end{theorem1}
\(A_p=1\) for \(2\leq p\) and \(B_p=1\) for \(p\leq 2\) follows directly from H\"older's inequality.
Haagerup published in \cite{Haagerup} the first complete proof  of Theorem \ref{MainTheorem}. A relatively elementary 
proof for the case \(p\geq3\) can be found in \cite{FHJSZ} where mainly convexity arguments were used. 
A new approach for the case \(0<p<2\) was given in \cite{NazarovPodkorytov} and we want to extend these ideas for the
case \(2<p<3\). 
\subsection{Organisation of the paper}
In section 2 we will recall the basic idea's to prove Khintchine's inequality which leads to the 
inequality of type \(\int(g^s-f^s)\mathrm{d}\mu\geq 0\). In section 3 we will describe the method of Nazarov
and Podkorytov to solve such kind of integral inequalities. Their lemma provides a positive answer to this inequality
by showing two conditions. Sections 4 and 5 will be dedicated to the proof of the first and second condition.
\subsection{Preliminaries}
Throughout this article, we only make use of standard notations and functions which should be well-known to
the reader.
At several points, we make use of the series expansion of the function \(f(t)=-\ln(\cos(t))\) for 
\(t\in(-\frac{\pi}{2},\frac{\pi}{2})\). We remark \(f'(t)=\tan(t)\) and the series representation of \(f\) is 
(cf. \cite{GradshteinRhyzhik}, 1.518):
\begin{align}
 -\ln(\cos(t))=\sum_{k=1}^\infty\frac{2^{2k-1}(2^{2k}-1)|B_{2k}|}{k(2k)!}t^{2k}
              =\frac{t^2}{2}+\frac{t^4}{12}+\frac{t^6}{45}+\dots\label{lncosRepresentation}
\end{align}
 where \(B_m\) denotes the m'th Bernoulli number. We remark that the coefficients of this
 series expansion are non-negative and hence, we can derive from this upper estimates for the cosine, i.e.:
\begin{align}
 \cos(t)\leq \exp(-t^2/2-t^4/12-t^6/45)\leq \exp(-t^2/2-t^4/12)\leq \exp(-t^2/2)\label{CosineInequality}
\end{align}
We will make use of the exponential, sine and cosine integral:
For \(x <0 \) we define
\begin{align}
 Ei(x)=-\int\limits_{-x}^\infty\frac{e^{-t}}{t}\mathrm{d}t=C+\ln(-x)+\sum_{k=1}^\infty\frac{x^k}{k\cdot k!}\label{EiFunction}
\end{align}
and for positive \(x\) we define
\begin{align}
 si(x)=-\int\limits_{x}^\infty \frac{\sin(t)}{t}\mathrm{d}t=-\frac{\pi}{2}-\sum_{k=1}^\infty\frac{(-1)^kx^{2k-1}}{(2k-1)(2k-1)!}\label{siFunction}
\end{align}
and
\begin{align}
 ci(x)=-\int\limits_x^\infty\frac{\cos(t)}{t}\mathrm{d}t=C+\ln(x)+\sum_{k=1}^\infty(-1)^k\frac{x^{2k}}{2k(2k)!}\label{ciFunction}
\end{align}
where \(C=0.577215\dots\) is Euler's constants (see \cite{GradshteinRhyzhik},8.211, 8.214 8.230, 8.232). 
We assume that one can always compute sufficiently 
precise numerical values of \(Ei(x)\), \(si(x)\) and \(ci(x)\) for fixed \(x\) via the series representation.

\section{Haagerup's approach}
The first observation
due to Ste\v{c}kin is that
\(\liminf\limits_{n\rightarrow\infty}\mathbbm{E}|\frac{1}{\sqrt{n}}\sum_{k=1}^n \varepsilon_k|^p 
\geq 2^{\frac{p}{2}}\frac{\Gamma(\frac{p+1}{2})}{\sqrt{\pi}}\) (see \cite{Steckin} and also \cite{Haagerup}). 
Hence, we only need to establish 
\(B_p\leq 2^{\frac{1}{2}}\left(\frac{\Gamma(\frac{p+1}{2})}{\sqrt{\pi}}\right)^{\frac{1}{p}}\). 
We introduce for \(s>0\) the following auxiliary function:
\begin{align}
 F_p(s)=c_p\int_0^\infty\left(\frac{1}{2}t^2-1+|\cos(t/\sqrt{s})|^s\right)\frac{\mathrm{d}t}{t^{p+1}}\notag
\end{align}
where \(c_p\) is a constant which is of no greater interest here.
We have the following lemma:
\begin{lemma1}\label{HaagerupLemma}
 \begin{enumerate}
  \item Let \(\sum_{k=1}^n a_k^2=1\) then we have
\begin{align}
 \mathbbm{E}|\sum_{k=1}^n a_k\varepsilon_k|^p\leq \sum_{k=1}^na_k^2F_p(a_k^{-2})
\end{align}
 \item 
\begin{align}
 \lim\limits_{s\rightarrow\infty}F_p(s)
=c_p\int_0^\infty\left(\frac{1}{2}t^2-1+\exp(-t^2/2)\right)\frac{\mathrm{d}t}{t^{p+1}}
=2^{\frac{p}{2}}\frac{\Gamma(\frac{p+1}{2})}{\sqrt{\pi}}\notag
\end{align} 
 \end{enumerate}
\end{lemma1}
If we were able to show 
\begin{align}
 F_p(s)\leq 2^{\frac{p}{2}}\frac{\Gamma(\frac{p+1}{2})}{\sqrt{\pi}}\label{inequality1}\qquad,
\end{align}
we would have with the help of the preceding lemma that for all sequences of scalars \(a_1, \dots, a_n\) with 
\(l_2\)-Norm 1:
\begin{align}
\mathbbm{E}|\sum_{k=1}^n a_k\varepsilon_k|^p\leq \sum_{k=1}^na_k^2F_p(a_k^{-2})
\leq 2^{p/2}\frac{\Gamma(\frac{p+1}{2})}{\sqrt{\pi}}\left(\sum_{k=1}^na_k^2\right)
=2^{p/2}\frac{\Gamma(\frac{p+1}{2})}{\sqrt{\pi}}\notag
\end{align}
and this establishes the Khintchine inequality by taking \(\sqrt[p]{\cdot}\) and by homogeneity.
Unfortunately, (\ref{inequality1}) is not true for \(s\) sufficiently close to \(1\) and \(p\) sufficiently close 
to \(2\). However, there is a short an elementary proof in \cite{Haagerup}, Lemma 4.6
for the case that \(\sum_{k=1}^n a_k^2 =1\) and that one of the coefficients \(|a_k|\geq 2^{\frac{1}{4}}\).  
So, we only need to consider the case that all
the coefficients are lesser or equal than \(2^{\frac{1}{4}}\), which means that we only need to show (\ref{inequality1}) for 
\(s\geq \sqrt{2}\) in order to show Khintchine's inequality. With Lemma \ref{HaagerupLemma}, 2. and the definition of 
\(F_p\) inequality (\ref{inequality1}) is equivalent to
\begin{align}
 \int\limits_0^\infty\left(\exp(-t^2/2)^s-|\cos(t)|^s\right)\frac{\mathrm{d}t}{t^{p+1}}\geq 0\label{inequality2}
\end{align}

\section{Nazarov's and Podkorytov's lemma on distribution function}
We introduce the definition of distribution functions suitable for our needs:
\begin{definition1}
 Let \(f:X\rightarrow \mathbbm{R}_{>0}\) a measurable function on a measure space \((X,\mu)\). Then we denote by
 \(F_*(y):=\mu(\{x\in X:f(x)<y\})\), \(y\in\mathbbm{R}_{>0}\),
the distribution function of \(f\)
\end{definition1}
Since we know that the integral of a function is determined by its distribution function, it might be more useful to 
study the distribution functions rather than the functions themselves in order to treat integral inequalities like
(\ref{inequality2}). In our case we can hope to simplify the situation because (\ref{inequality2}) has an oscillating
integrand whereas the distribution functions are increasing. In \cite{NazarovPodkorytov} we find a lemma
on distribution functions which was successful used in this article to compute the optimal constants for the case \(0<p<2\).
Later, it was also used in \cite{Koenig} for the optimal constants in the Khintchine inequality with Steinhaus variables.
It states the following:

\begin{lemma1}\label{DistributionFunction}
 Let \(Y>0\) and \(f,g :X\rightarrow [0,Y]\) measurable functions on a measure space \((X,\mu)\).
We assume the distribution functions \(F_*,G_*\) of \(f,g\) to be finite on \((0,Y)\). Let 
\(S=\{s>0:g^s-f^s\in L_1(X,\mu)\}\). If there exist \(y_0\in(0,Y)\) and \(s_0\in S\) such that the following two conditions
hold:
\begin{enumerate}
 \item \(F_*-G_*\leq 0\) on \((0,y_0)\)
       and \(F_*-G_*\geq 0\) on \((y_0,Y)\)
 \item \(\int\limits_X\left(g^{s_0}-f^{s_0}\right)\mathrm{d}\mu\geq 0\qquad,\)
\end{enumerate}
then we have
\(\int\limits_X\left(g^s-f^s\right)\mathrm{d}\mu\geq 0\)
for all \(s\in S\) with \(s\geq s_0\).
\end{lemma1}

In our situation we will consider the measure space
\(X=(0,\infty)\) with the measure \(\mathrm{d}\mu_p=\frac{\mathrm{d}t}{t^{p+1}}\) and the
two functions \(f,g\) on \(X\) which are defined by \(f(t):=|\cos(t)|\) and \(g(t)=\exp(-t^2/2)\).
An easy calculation shows that the distribution functions of \(f, g\) are:
\begin{align}
 F_*(x)=\frac{1}{p}\sum\limits_{k=0}^\infty\left(\frac{1}{(k\pi+\arccos(x))^p}-\frac{1}{((k+1)\pi-\arccos(x))^p}\right)\label{FormulaF}
\end{align}
and
\(G_*(x)=\frac{1}{p}(-2\ln(x))^{-p/2}\).
We take \(Y=1\) and one can verify that in terms of Lemma \ref{DistributionFunction} we have \(S=\mathbbm{R}_{>0}\) and
\(F_*, G_*\) are finite functions and \(C^1\) on \((0,1)\).

\section{Proof of the first condition of lemma \ref{DistributionFunction}}
We will proof exactly the following lemma from which the first condition follows:
\begin{lemma1}\label{Hauptlemma}
 There exists constants \(0<\rho<\sigma<1\) such that the following three properties are fulfilled:
\begin{enumerate}
 \item \(F_*(\sigma)-G_*(\sigma)\geq0\)
 \item \(F_*-G_*\) is increasing on \((\rho,1)\)
 \item \(F_*-G_*\) is lesser than \(0\) on \((0,\rho)\)
\end{enumerate}
\end{lemma1}
We will show this lemma with \(\rho=\frac{1}{15}\) and \(\sigma=0.97\).

\subsection{Proof of Lemma \ref{Hauptlemma}, 1.}  
We minorize \(F_*\):
\begin{align}
 pF_*(x)&=\frac{1}{\arccos(x)^p}-\sum\limits_{k\in\mathbbm{N}}\left(\frac{1}{(k\pi-\arccos(x))^p}-
                                                              \frac{1}{(k\pi+\arccos(x))^p}\right)\notag\\
         &\geq \frac{1}{\arccos(x)^p}-\sum\limits_{k\in\mathbbm{N}}\left(\frac{1}{(k\pi-k\arccos(x))^p}-
                                                              \frac{1}{(k\pi+k\arccos(x))^p}\right)\notag\\
         &=\frac{1}{\arccos(x)^p}-\left(\sum\limits_{k\in\mathbbm{N}}\frac{1}{k^p}\right)\left(\frac{1}{(\pi-\arccos(x))^p}-
                                                              \frac{1}{(\pi+\arccos(x))^p}\right)\notag\\
         &\geq\frac{1}{\arccos(x)^p}-\left(\sum\limits_{k\in\mathbbm{N}}\frac{1}{k^2}\right)\left(\frac{1}{(\pi-\arccos(x))^p}-
                                                              \frac{1}{(\pi+\arccos(x))^p}\right)\label{LastLine1}
\end{align}
And the series in (\ref{LastLine1}) equals \(\frac{\pi^2}{6}\). By the mean value theorem applied to the function 
\(x\mapsto-\frac{1}{x^p}\), there is a \(\zeta\in(\pi-\arccos(x), \pi+\arccos(x))\) such that:
\begin{align}
 &\frac{1}{2\arccos(x)}\left(\frac{1}{(\pi-\arccos(x))^p}-\frac{1}{(\pi+\arccos(x))^p}\right)\notag\\
=&p\frac{1}{\zeta^{p+1}}\leq 3\frac{1}{(\pi-\arccos(x))^{p+1}}
\leq 3\frac{1}{(\pi-\arccos(x))^3}\notag
\end{align}
Hence, (\ref{LastLine1}) is greater or equal than
\begin{align}
 \frac{1}{\arccos(x)^p}-\pi^2\frac{\arccos(x)}{\left(\pi-\arccos(x)\right)^3}\notag
\end{align}
and therefore, we get:
\begin{align}
 p(F_*(x)-G_*(x))\geq \frac{1}{\arccos(x)^p}-\frac{1}{\sqrt{2\cdot\ln(\frac{1}{x})}^p}-\pi^2\frac{\arccos(x)}{(\pi-\arccos(x))^3}\label{Inequality2}
\end{align}
If we substitute with \(x=0.97\) and \(p=2\), we see that the right-hand side is greater than \(0\). Since
\(\frac{1}{\arccos(0.97)}>\frac{1}{\sqrt{2\cdot\ln(\frac{1}{0.97})}}\geq 1\), the right-hand side of (\ref{Inequality2}) 
increases for \(x=0.97\),
if we increase \(p\) and this shows 3.
\hfill\(\Box\)

\subsection{Proof of Lemma \ref{Hauptlemma}, 2.} We recall \(2<p<3\) and \(0<x< \frac{1}{15}\). We can represent the summands of \(F_*\) in the 
following manner:
\begin{align}
 \frac{1}{(k\pi+\arccos(x))^p}-\frac{1}{((k+1)\pi-\arccos(x))^p}
=\frac{\frac{1}{(1-\varepsilon_k)^p}-\frac{1}{(1+\varepsilon_k)^p}}{[(k+\frac{1}{2})\pi]^p}\notag
\end{align}
 where \(\varepsilon_k=\frac{\pi/2-\arccos(x)}{(k+1/2)\pi}\leq 0.04248\). Using linear Taylor approximation and 
estimating the remainder one gets
\begin{align}
 (1+\varepsilon_k)^p-(1-\varepsilon_k)^p\leq 2p(1+\varepsilon_k^2)\varepsilon_k\leq2.00361p\varepsilon_k\notag
\end{align}
We remark that \(\pi/2-\arccos(x)\leq 0.06672\) and \(\sin(\pi/2-\arccos(x))=x\). By concavity of the sine on the 
interval \((0,\pi)\) we get for \(0\leq t\leq 0.06672\) that
\(\sin(t)\geq \frac{t}{0.06672}\sin(0.06672)+\left(1-\frac{t}{0.06772}\right)\sin(0)\) and thus,
\(t\leq \frac{0.06672}{\sin(0.06672)}\sin(t)\notag
\).
With all these results in mind we calculate:
\begin{align}
&\frac{1}{(1-\varepsilon_k)^p}-\frac{1}{(1+\varepsilon_k)^p}
\leq \frac{(1+\varepsilon_k)^p-(1-\varepsilon_k)^p}{(1-\varepsilon_k^2)^p}
\leq \frac{(1+\varepsilon_k)^p-(1-\varepsilon_k)^p}{(1-0.04248^2)^3}\notag\\
\leq&2.0145 p \varepsilon_k
\leq 2.0145p\cdot\frac{0.06672}{\sin(0.06772)}\frac{\sin(\pi/2-\arccos(x))}{(k+\frac{1}{2})\pi}
\leq \frac{2.02}{(k+\frac{1}{2})\pi}px\label{LetzteZeile3}
\end{align}
Replacing the expression (\ref{LetzteZeile3}) in (\ref{FormulaF}), we get:
\begin{align}
 F_*(x)\leq 2.02\left(\sum\limits_{k=0}^\infty\frac{1}{[(k+\frac{1}{2})\pi]^{p+1}}\right)\cdot x=:d_p\cdot x\notag
\end{align}
The coefficient \(d_p\) is convex in \(p\), which yields \(d_p\leq d_2+(d_3-d_2)(p-2)\). We refind via the 
representation
\begin{align}
 d_p=2.02\left(\frac{2}{\pi}\right)^{p+1}\left(1-\frac{1}{2^{p+1}}\right)\left(\sum\limits_{k=1}^\infty\frac{1}{k^{p+1}}\right)\notag
\end{align}
Riemann's \(\zeta\)-function, so we can compute with sufficient precision \(d_2\leq0.5482\) and \(d_3\leq0.3367\)
and we conclude
\(F_*(x)\leq(0.98-0.2115p)\cdot x\).
The last step will be to show \((0.98-0.2115p)\cdot x\leq G_*(x)\) which is equivalent to
\(p(0.98-0.2115p)\leq\frac{x^{-1}}{(2\ln(1/x))^{p/2}}\)
or 
\begin{align}
 p(0.98-0.2116p)\leq \frac{\exp(t)}{(2t)^{p/2}}\label{Inequality3}
\end{align}
for \(t\geq 2.7\).
Maximizing the left-hand side of (\ref{Inequality3}) in \(2\leq p \leq 3\) we get \(p(0.98-0.2115p)\leq1.14\).
The right-hand side of (\ref{Inequality3}) is increasing on \([\frac{3}{2},\infty)\) 
(which can be verified via the derivation), so we have
\(\frac{\exp(t)}{(2t)^{p/2}}\geq\frac{\exp(2.7)}{(2\cdot 2.7)^{3/2}}\geq 1.8\geq 1.4\).
\hfill\(\Box\)

\subsection{Proof of Lemma, \ref{Hauptlemma}, 3.} We prove that \(F_*-G_*\) is monotonously 
 increasing on \((\frac{1}{15},1]\) which is equivalent to \(\frac{F'_*}{G'_*}\geq 1\) providing \(G'_*>0\). We 
compute
\begin{align}
 F'_*(x)&=\sum\limits_{k=0}^\infty\left(\frac{1}{(k\pi+\arccos(x))^{p+1}}+\frac{1}{((k+1)\pi-\arccos(x))^{p+1}}\right)
          \frac{1}{\sqrt{1-x^2}}\notag\\
        &\geq \left(\frac{1}{\arccos(x)^{p+1}}+\frac{1}{(\pi-\arccos(x))^{p+1}}\right)\frac{1}{\sqrt{1-x^2}}\notag
\end{align}
and
\(G'_*(x)=\frac{1}{x(-2\ln(x))^{p/2+1}}\notag\).
For simplification, we set \(t=\arccos(x)\). This yields
\begin{align}
 \frac{F'_*(x)}{G'_*(x)}&\geq \left(\frac{1}{t^{p+1}}+\frac{1}{(\pi-t)^{p+1}}\right)\left(-2\ln(\cos(t))\right)^{\frac{p}{2}+1}
                              \cot(t)\notag\\
&= \left(\sqrt{\frac{-2\ln(\cos(t))}{t^2}}^{p+1}+\sqrt{\frac{-2\ln(\cos(t))}{(\pi-t)^2}}^{p+1}\right)
                  (-2\ln(\cos(t)))^{\frac{1}{2}}\cot(t)\label{LetzteZeile4}
\end{align}
Hence, it suffices to show that the right-hand side is greater than 1 for \(0< t < 1.50409\).
We split up our proof in subchapters.

\subsubsection{Reduction to the case $p=2$}
We show that the case \(p=2\) gives the least value in (\ref{LetzteZeile4}) for every \(t\) and for this, it 
is appropriate to examine the expression \(A^{p+1}+B^{p+1}\). We have the following lemma which easy to check:
\begin{lemma1}
 Let \(p\geq 2\), \(A\geq 1\) and \(B>0\) be real numbers. If the inequality \(\frac{A^3}{B^3}\cdot\ln(A)\geq -\ln(B)\) 
holds, then we have \(A^{p+1}+B^{p+1}\geq A^3+B^3\).
\end{lemma1}
Remark that the case \(B\geq 1\) is trivial.
We put \(A:=\sqrt{\frac{-2\ln(\cos(t))}{t^2}}\) and \(B:=\sqrt{\frac{-2\ln(\cos(t))}{(\pi-t)^2}}\).
\(A \geq 1\) follows from the fact that \(-\ln(\cos(t))\geq \frac{t^2}{2}\) on \([0, \frac{\pi}{2}]\) 
(see (\ref{lncosRepresentation})). We have \(\frac{A^3}{B^3}=\frac{(\pi-t)^3}{t^3}\)
and hence, we want to show the inequality
\begin{align}
 (\pi-t)^3\ln\left(\frac{-2\ln(\cos(t))}{t^2}\right)\geq -t^3\ln\left((\frac{-2\ln(\cos(t))}{(\pi-t)^2}\right)\label{LetzteZeile5}
\end{align}

Bringing the terms, which belong to \(t^3\), on the right-hand side and applying the logarithm rules, we  can restate
(\ref{LetzteZeile5}) as follows:
\begin{align}
 (\pi^3-3\pi^2t+3\pi t^2)\ln\left(\frac{-2\ln(\cos(t))}{t^2}\right)\geq 2t^3\ln\left(\frac{\pi-t}{t}\right)\notag
\end{align}
Use \(\frac{-2\ln(\cos(t))}{t^2}\geq 1+\left(\frac{t^2}{6}+\frac{2}{45}t^4\right)\)
for \(0\leq t < \frac{\pi}{2}\) (see \ref{lncosRepresentation}) and the well-known estimate  \(\ln(1+x)\geq x-\frac{x^2}{2}\) for \(x\geq 0\) and conclude:
\begin{align}
 \ln\left(\frac{-2\ln(\cos(t))}{t^2}\right)
\geq \frac{t^2}{6}+t^4\left(\frac{2}{45}-\frac{1}{2}\left[\frac{1}{6}+\frac{2}{45}t^2\right]^2\right)\label{LetzteZeile6}
\end{align}
Since \(t\) is bounded from above by \(\frac{\pi}{2}\), we can check that the coefficient of \(t^4\) is
positive, which yields that the left-hand side of (\ref{LetzteZeile6}) is bounded from below
by \(\frac{t^2}{6}\). Now, it is sufficient to proof that 
\((\pi^3-3\pi^2t+3\pi t^2)\cdot \frac{t^2}{6}\geq 2t^3\ln(\frac{\pi-t}{t})\) or equivalently, 
\begin{align}
\pi^3-3\pi^2t+3\pi t^2\geq 12 t \ln\left(\frac{\pi-t}{t}\right)\label{LetzteZeile7}
\end{align} 
We will describe how one can show this inequality. 
First, we remark that \(t\ln(\frac{\pi-t}{t})\) is concave which can be proved by looking at the 
derivative, which is monotonously decreasing. Denote by \(T_{t_0}\) the function of the tangent
to \(12t\ln(\frac{\pi-t}{t})\) at the point \(t_0\). By concavity, we have
\begin{align}
 12t\ln\left(\frac{\pi-t}{t}\right)\leq T_{t_0}(t)
\end{align}
for every \(t_0\). Replacing this estimation in (\ref{LetzteZeile7}) leads to the inequality
\(\pi^3-3\pi^2t+3\pi t^2\geq T_{t_0}(t)\) which is obviously an inequality  with only second degree polynomials.
Taking \(t_0=1\) this inequality can be verified.
\hfill\(\Box\)
\\\\
Thus, we have shown that it suffices to proof
\begin{align}
 \left(\frac{1}{t^3}+\frac{1}{(\pi-t)^3}\right)[-2\ln(\cos(t))]^2\cot(t)\geq 1\label{LetzteZeile8}
\end{align}
and we will treat the case \(0\leq t \leq 1\) and \(t>1\) separately and by different methods
\subsubsection{Case 1: $0\leq t \leq 1$ and $p=2$}
We will minorize every of the three factors in (\ref{LetzteZeile8}) by a polynomials in \(t\) and in \(\frac{1}{t}\), 
which is described exactly by the following lemma:
\begin{lemma1} Let \(0<t\leq 1\) be a recall number, then the three following inequalities hold:
 \begin{enumerate} 
  \item 
  \(\frac{1}{t^3}+\frac{1}{(\pi-t)^3}\geq \frac{1}{t^3}\left(1+\frac{1}{\pi^3}t^3+\frac{3}{\pi^4}t^4+\frac{6}{\pi^5}t^6\right)\)
  \item 
  \([-2\ln(\cos(t))]^2\geq t^4\left(1+\frac{1}{3}t^2+\frac{7}{60}t^4\right)\)
  \item \(\cot(t)\geq \frac{1}{t}-\frac{t}{3}-\frac{1}{40}t^3\) 
\end{enumerate}
\end{lemma1}
 \underline{Proof of 1.:} We can transform the left-hand side to 
\(\frac{1}{t^3}\left(1+\frac{t^3}{\pi^3}\left(\frac{1}{1-t/\pi}\right)^3\right)\). Expanding \(\frac{1}{1-t/\pi}\)
to a geometric series, we minorize as follows:
\begin{align}
 \left(\frac{1}{1-t/\pi}\right)^3
=\left(\sum\limits_{k=0}^\infty\left(\frac{t}{\pi}\right)^k\right)^3
\geq\left(1+\frac{1}{\pi}t+\frac{1}{\pi^2}t^2\right)^3
\geq 1+\frac{3}{\pi}t+\frac{6}{\pi^2}t^2\notag
\end{align}
which gives the desired result.
\\\\
\underline{Proof of 2.:} We use \(\frac{-2\ln(\cos(t))}{t^2}\geq 1+\frac{t^2}{6}+\frac{2}{45}t^4\) 
(see (\ref{lncosRepresentation})) and we evaluate
\((1+\frac{t^2}{6}+\frac{2}{45}t^4)^2\) up to the power \(4\). 
\\\\
\underline{Proof of 3.:} The series expansion of the cotangent  (see \cite{GradshteinRhyzhik}) is:
\begin{align}
 \cot(t)=\frac{1}{t}-\frac{1}{3}t
-t^3\underbrace{\sum\limits_{k=2}^\infty\frac{2^{2k}|B_{2k}|}{(2k)!}t^{2k-4}}_{=:R(t)}\notag
\end{align}
Since \(R\) is increasing on \([0,1]\), we have \(R(t)\leq R(1)=\frac{1}{1}-\frac{1}{3}-\cot(1)\leq \frac{1}{40}\).
\hfill\(\Box\)
\begin{proposition1}
 Let \(0\leq t \leq 1\), then we have:
\begin{align}
 \left(1+\frac{1}{\pi^3}t^3+\frac{3}{\pi^4}t^4+\frac{6}{\pi^5}t^5\right)\left(1-\frac{1}{3}t^2-\frac{1}{40}t^4\right)
\geq 1-\frac{1}{3}t^2+\frac{1}{40}t^3\label{LetzteZeile11}
\end{align}
\end{proposition1}
\underline{Proof:} We expand partially the left-hand side of (\ref{LetzteZeile11}) such that we can represent this as
the sum of \(p_1(t):=1-\frac{1}{3}t^2+\frac{1}{\pi^3}t^3+\frac{3}{\pi^4}t^4\) and
\begin{align}
 p_2(t):=t^5\left(\frac{6}{\pi^5}-\left[\frac{1}{3\pi^3}+\frac{1}{\pi^4}t+(\frac{2}{\pi^5}+\frac{1}{40\pi^3})t^2
                                +\frac{3}{40\pi^4}t^3+\frac{3}{20\pi^5}t^4\right]\right)\label{LetzteZeile12}
\end{align}
We can maximize the expression in square brackets by simply taking \(t=1\) and thus, we can bound (\ref{LetzteZeile12})
from below by \(-0.009t^5\geq -0.009\cdot 1\cdot t^4\). Now, we have 
\(p_1(t)+p_2(t)\geq p_1(t)-0.009t^4\geq 1-\frac{1}{3}t^2+\frac{1}{\pi^3}t^3-0.004t^4\) and we show analogously as above the estimation
\(-0.004t^4\geq -0.004t^3\) from which we obtain the desired inequality (\ref{LetzteZeile11}).
\hfill\(\Box\)
\begin{corollary1}
Let \(0\leq t\leq 1\), then we have:
\begin{align}
 \left(1-\frac{1}{3}t^2+\frac{1}{40}t^3\right)\left(1+\frac{1}{3}t^2+\frac{7}{60}t^4\right)\geq 1\label{LetzteZeile13}
\end{align}
\end{corollary1}
\underline{Proof:} The expansion of the left-hand side of (\ref{LetzteZeile13}) is exactly
\(1+\frac{11}{360}t^4-\frac{11}{360}t^6+\frac{7}{2400}t^8\). This expression is obviously
greater or equal to \(1\) for \(0\leq t \leq 1\).
\hfill\(\Box\)
\subsubsection{Case 2: $1\leq t \leq 1.50409$ and $p=2$} The desired
inequality can be reformulated as 
\begin{align}
 g(t):=\frac{1}{t^3}+\frac{1}{(\pi-t)^3}\geq \frac{\tan(t)}{[-2\ln(\cos(t))]^2}=:f(t)\notag
\end{align}
\(g\) is obviously a convex function on \((0,\frac{\pi}{2})\) and one can show that \(f\) is convex there, too. 
Denote by \(T_{t_0}\) the function of tangent to \(g\) at the point \(t_0\). We can check that \(g\geq f\) on an interval \([a,b]\)
by the following method:
\\
Find an appropriate \(t_0\) and show \(T_{t_0}(a)\geq f(a)\) and as well for \(b\). Thus, we have
by convexity \(g(t)\geq T_{t_0}(t)\geq f(t)\) for all \(a\leq t \leq b\).
\\
Good values to do so are \(t_0=1.1\) for the interval \([1, 1.25]\) and \(t_0=1.45\) for the interval \([1.24,1.505]\).
\\[1ex]
We now give a sketch of the proof that \(f\) is convex on \((0,\frac{\pi}{2})\):
\\[2ex]
\underline{Proof:} We show that \(f''(t)\geq 0\). The second derivative of \(f''(t)\) multiplied
by the positive value \(\frac{2\ln(\cos(t))^4\cos(t)^3}{\sin(t)}\)
is \(\ln(\cos(t))^2+3\ln(\cos(t))+3\sin(t)^2\). We put \(s:=-\ln(\cos(t))\) and we conclude that 
\(f''(t)\geq 0\) for all \(0<t<\pi/2\) is equivalent to \(s^2-3s+3-3\exp(-2s)\geq 0\) for all \(0<s<\infty\).
We multiply both sides by \(\exp(2s)\) and we use \(\exp(2s)\geq 1+2s+2s^2\). Then, we have
\begin{align}
 (s^2-3s+3)(1+2s+2s^2)-3=s(2s(s-1)^2+3-s)\notag
\end{align}
which is indeed positive for positive \(s\). 
\hfill\(\Box\)

\section{Proof of the second condition of lemma \ref{DistributionFunction}}
We need to show that 
\begin{align}
 H(p):=\int\limits_{0}^\infty \frac{\exp(-t^2/\sqrt{2})-|\cos(t)|^{\sqrt{2}}}{t^{p+1}}\mathrm{d}t\geq 0\notag
\end{align}
for all \(2<p<3\). We do this by showing the two 
properties \(H(2)\geq 0\) and \(H'(p)\geq 0\) for \(2\leq p \leq 3\). 
\subsection{Proof of $H'(p)\geq 0$}
We have by differentiation under the integral
\begin{align}   
 H'(p)=\int_0^\infty -\ln(t)\frac{\exp(-t^2/\sqrt{2})-|\cos(t)|^{\sqrt{2}}}{t^{p+1}}\mathrm{d}t\notag
\end{align}
and we examine the integral on different intervals. Let us start with \([0,1]\). Since \(-\ln(t)\geq 0\) and
\(\exp(-\frac{t^2}{\sqrt{2}})\geq \exp(-\frac{t^2}{\sqrt{2}}-\frac{t^4}{\sqrt{2}\cdot6})\geq |\cos(t)|^{\sqrt{2}}\) and
we obtain the estimate
\begin{align}
 -\ln(t)\frac{\exp(-t^2/\sqrt{2})-|\cos(t)|^{\sqrt{2}}}{t^{p+1}}\geq -\ln(t)\frac{\exp(-t^2/\sqrt{2})
                                                                  \left(1-\exp(-\frac{t^4}{\sqrt{2}\cdot 6})\right)}{t^3}\label{LetzteZeile20}
\end{align}
We estimate each of the three (positive) factors on the right-hand side of (\ref{LetzteZeile20}) from below by Taylor polynomials.
We use \(-\ln(t)\geq 1-t\) and additionally \(\exp(-a)\geq 1-a\) and \(1-\exp(-b)\geq b-\frac{b^2}{2}\) for \(a,b\geq 0\) which yield sufficiently
precise estimates for the second and third factor. The result of these estimates for the
integrand is a not to extensive polynomial, for which we can compute the exact value \(0.0153\).
\\\\
The next interval is \([1,\frac{\pi}{2}]\). The integrand is negative since \(-\ln(t)\leq 0\) and
\(\exp(-\frac{t^2}{\sqrt{2}})\geq |\cos(t)|^{\sqrt{2}}\) and consequently, we want to find a precise estimate from above
of
\begin{align}
 \frac{\ln(t)\left(\exp(-\frac{t^2}{\sqrt{2}})-|\cos(t)|^{\sqrt{2}}\right)}{t^{p+1}}
 \leq \frac{\ln(t)\left(\exp(-\frac{t^2}{\sqrt{2}})-|\cos(t)|^{\sqrt{2}}\right)}{t^{3}}\notag
\end{align}
for \(1\leq t\leq \frac{\pi}{2}\).
We make use of the fact that \(\ln(t)\leq t-1\). Since \(f(t):=\exp(-t^2/\sqrt{2})\) is convex on this interval,
we can estimate this term from above by the secant function through \((1,f(1))\) and \((\frac{\pi}{2},f(\frac{\pi}{2}))\).
We estimate \(|\cos(t)|^{\sqrt{2}}\) from below by \(\cos(1.2)^{\sqrt{2}}\) on the interval \([1, 1.2]\) and 
by \(\cos(1.4)^{\sqrt{2}}\)
on the interval\([1.2, 1.4]\) and \(0\) elsewhere. These approximations give
\begin{align}
 \int\limits_1^{\frac{\pi}{2}}-\ln(t)\frac{\exp(-t^2/\sqrt{1})-|\cos(t)|^{\sqrt{2}}}{t^{p+1}}\mathrm{d}t\geq -0.0147\notag
\end{align}
It should be natural that the contribution to the integral on the interval \([\frac{\pi}{2},\infty)\)
should be positive, because \(\exp(-\frac{t^2}{\sqrt{2}})\) decays very fast whereas \(|\cos(t)|^{\sqrt{2}}\)
oscillates between \(0\) and \(1\).
By simple one-variable analysis methods one can show the following lemma
\begin{lemma1} Let \(t\in [\frac{\pi}{2},\infty)\), \(2\leq p \leq 3\), set \(\lambda_p:=1.75\cdot\left(\frac{2}{\pi}\right)^p\)
and \(\Lambda_p=\frac{1}{e(p-1)}\), then we have
 \begin{enumerate}
  \item \(\frac{\ln(t)}{t^{p+1}}\geq \lambda_p\frac{1}{t^4}\)
  \item \(\frac{\ln(t)}{t^{p+1}}\leq \Lambda_p \frac{1}{t^2}\)
 \end{enumerate}
\end{lemma1}
The motivation behind this lemma is that we are now able to give estimates which we can compute explicitly 
(see \cite{GradshteinRhyzhik}, 2.642):
\begin{align}
 \int\limits_{\frac{\pi}{2}}^{\infty}\ln(t)\frac{|\cos(t)|^{\sqrt{2}}}{t^{p+1}}
\geq \lambda\int\limits_{\frac{\pi}{2}}^\infty \frac{\cos(t)^2}{t^4}\mathrm{d}t\notag
= \lambda_p \frac{2+\pi^2-2si(\pi)\pi}{3\pi}\geq 0.043369 \left(\frac{2}{\pi}\right)^p\notag
\end{align}
where \(si\) was defined in (\ref{siFunction}).
On the other hand, we have:
\begin{align}
 \int\limits_{\frac{\pi}{2}}^\infty \ln(t)\frac{\exp(-t^2/\sqrt{2})}{t^{p+1}}\mathrm{d}t
\leq \Lambda_p \int_{\frac{\pi}{2}}^\infty\frac{\exp(-t^2/\sqrt{2})}{t^2}\mathrm{d}t
\leq 0.00705 \frac{1}{p-1}\notag
\end{align}
The  second integral was computable since we can transform it via partial integration to an integral
with integrand of the form \(\exp(-s^2)\).
Hence, we need to show that \(0.043369\left(\frac{2}{\pi}\right)^p\geq 0.00705 \frac{1}{p-1}\).
This is feasible because \((p-1)\left(\frac{2}{\pi}\right)^p\) is increasing on \([2,3]\) and so one only needs to check this 
for \(p=2\). This shows that the integral on the interval \([\frac{\pi}{2},\infty)\) is positive.

\subsection{Proof of $H(2)\geq 0$}
We want to show that 
\begin{align}
 \int\limits_0^\infty \frac{\exp(-t^2/\sqrt{2})-|\cos(t)|^{\sqrt{2}}}{t^3}\mathrm{d}t\geq 0\notag
\end{align}
and we start with  our estimates on the interval \([0,\frac{\pi}{4}]\). We have 
\(|\cos(t)|^{\sqrt{2}}\leq \exp(-\frac{t^2}{\sqrt{2}}-\frac{t^4}{\sqrt{2}\cdot 6}-\frac{\sqrt{2}t^6}{45})\) by (\ref{CosineInequality})
and \(1-\exp(-b)\leq b-\frac{b^2}{2}\) for \(a,b\geq 0\) which yields
\begin{align}
 &\exp(-t^2/\sqrt{2})-|\cos(t)|^{\sqrt{2}}
\geq \exp(-t^2/\sqrt{2})\left(1-\exp\left(-\left[\frac{t^4}{\sqrt{2}\cdot 6}+\frac{\sqrt{2}t^6}{45}\right]\right)\right)\notag\\
\geq &\exp(-t^2/\sqrt{2})\left(\left[\frac{1}{\sqrt{2}6}t^4+\frac{\sqrt{2}}{45}t^6\right]-\frac{1}{2}\left[\frac{1}{\sqrt{2}\cdot 6}t^4+\frac{\sqrt{2}}{45}t^6\right]^2\right)\notag\\
\geq &\exp(-t^2/\sqrt{2})\left(\left[\frac{1}{\sqrt{2}\cdot 6}t^4+\frac{\sqrt{2}}{45}t^6\right]
-\frac{1}{2}\left(\frac{\pi}{4}\right)^2\left[\frac{\sqrt{2}}{12}+\frac{\sqrt{2}}{45}\left(\frac{\pi}{4}\right)^2\right]^2t^6\right)\label{LetzteZeile21}
\end{align}
Dividing the expression in (\ref{LetzteZeile21}) by \(t^3\), we get an expression of the form 
\(\exp(-\frac{t^2}{\sqrt{2}})(at+bt^3)\). Applying the transform of coordinates \(s\mapsto \sqrt[4]{2}\sqrt{s}\)
we end up with an expression of the form \(\exp(-t)(a'+b't)\) for which we can easily calculate the integral and
which gives finally the estimate:
\begin{align}
 \int\limits_0^{\frac{\pi}{4}} \frac{\exp(-t^2/\sqrt{2})-|\cos(t)|^{\sqrt{2}}}{t^3}\mathrm{d}t\geq 0.03129\notag
\end{align}
Next, we want to determine the exact integral of \(\exp(-t^2/\sqrt{2})/t^3\). For this, let
\(a>0\). Applying \(s\mapsto \sqrt[4]{2}\sqrt{s}\) and \cite{GradshteinRhyzhik}, 2.325, we get:
\begin{align}
 \int\limits_{a}^\infty \frac{e^{-t^2/\sqrt{2}}}{t^3}\mathrm{d}t
=\frac{1}{2\sqrt{2}}\int\limits_{\frac{a^2}{\sqrt{2}}}^\infty\frac{e^{-s}}{s^2}\mathrm{d}s
=\frac{e^{-a^2/\sqrt{2}}}{2a^2}+\frac{Ei\left(-\frac{a^2}{\sqrt{2}}\right)}{2\sqrt{2}}\label{LetzteZeile22}
\end{align}
where \(Ei\) is defined in (\ref{EiFunction}).
Substituting \(a=\frac{\pi}{4}\) we get, that (\ref{LetzteZeile22}) is bigger than \(0.29587\).
\\[2ex]
The last part is to estimate
\(\int\limits_{\frac{\pi}{4}}^\infty \frac{|\cos(t)|^{\sqrt{2}}}{t^3}\)
and we treat the intervals \([\frac{\pi}{4},\frac{3}{4}\pi]\)  and \([\frac{3}{4}\pi,\infty)\) separately.
Let us start with the interval \([\frac{\pi}{4}, \frac{3}{4}\pi]\). We will make use of the fact that
\begin{align}
 \int\frac{\cos(t)}{t^3}\mathrm{d}t= \frac{1}{2}\left(-\frac{\cos(t)}{t^2}+\frac{\sin(t)}{t}-ci(t)\right)+K\notag
\end{align}
and 
\begin{align}
 \int \frac{\cos(t)^2}{t^3}\mathrm{d}t= \frac{1}{4}\left(\frac{\cos(2t)}{t^2}+2\frac{\sin(2t)}{t}-4ci(2t)-\frac{1}{t^2}\right)\label{LetzteZeile23}
\end{align}
(see for example \cite{GradshteinRhyzhik}, 2.642), where \(ci\) was defined in (\ref{ciFunction}).
We have \(|\cos(t)|\in[0, \frac{\sqrt{2}}{2}]\) for \(t\in [\frac{\pi}{4}, \frac{3}{4}\pi]\) and our strategy will be two
bound \(x^{\sqrt{2}}\) from above with sufficient precision by  a polynomial  of the form \(a+bx+cx^2\) since we know 
how to integrate \(1/t^3\), \(\cos(t)/t^3\) and \(\cos(t)^2/t^3\).
\begin{lemma1}
 Let \(x\) be real and non-negative. We have:
\begin{enumerate}
 \item \(0 \leq (\sqrt{2}-1)x^2+(2-\sqrt{2}-0.126)x-x^{\sqrt{2}}\) for \(x\in[0,0.25]\)
 \item \(0 \leq (\sqrt{2}-1)x^2+0.6355x-0.04399-x^{\sqrt{2}}\) for \(x\in[0.25,\frac{\sqrt{2}}{2}]\).
\end{enumerate}
\end{lemma1}
\underline{Proof of 1.:} The function \(f(x)=(\sqrt{2}-1)x^2+(2-\sqrt{2})x-x^{\sqrt{2}}\) is concave on \([0,0.25]\)
and \(f(0)=0<f(0.25)\), which yields \(f(x)\geq0\) for all \(x\in[0,0.25]\). Hence, we can subtract
the function of the secant through the points \((0,f(0))\) and \((0.25, f(0.25))\) from \(f\)  and this gives the
nonnegative function \(x\mapsto\sqrt{2}-1)x^2+(2-\sqrt{2}-0.126)x-x^{\sqrt{2}}\).
\\[1ex]
\underline{Proof of 2.:} Set \(g(x):=(\sqrt{2}-1)x^2+0.6355x-0.04399-x^{\sqrt{2}}\) for 
\(x\in[0.25,\frac{\sqrt{2}}{2}]\). In order to show, that \(g\geq 0\) we remark \(g(\frac{\sqrt{2}}{2})\geq0\)
and our claim will follow from \(0\geq g'(x)=2(\sqrt{2}-1)x+0.6335-\sqrt{2}x^{\sqrt{2}-1}\).
Since \(g'\) is convex, we are finished by the fact that \(g'(0.25)\leq0\) and \(g'(\frac{\sqrt{2}}{2})\leq0\).
\hfill\(\Box\)
\\[2ex]
Using this lemma to bound \(|\cos(t)|^{\sqrt{2}}\) from above by a polynomial of the form
\(a\cos(t)^2+ b|\cos(t)|+c\) we get the following estimate: 
\(\int_{\frac{1}{4}\pi}^{\frac{3}{4}\pi}\frac{|\cos(t)|^{\sqrt{2}}}{t^3}\mathrm{d}t\leq 0.2577\).
Now, we want to bound \(\int_{\frac{3}{4}\pi}^\infty \frac{|\cos(t)|^{\sqrt{2}}}{t^3}\mathrm{d}t\)
from above. For this, we consider the measure space \(X=(\frac{3}{4}\pi, \infty)\) endowed with the measure
\(\mathrm{d}\mu:=\frac{\mathrm{d}t}{t^3}\). Since \(X\) has finite measure, we can apply the following 
well-known H\"older-type inequality:
\begin{align}
 \|\cos(t)\|_{L_{\sqrt{2}}(\mu)}\leq \mu(X)^{\frac{1}{\sqrt{2}}-\frac{1}{2}}\|\cos(t)\|_{L_2(\mu)}\notag
\end{align}
Using the primitive of \(\frac{\cos(t)^2}{t^3}\) (see \ref{LetzteZeile23}), we can compute the 
right-hand side. From this, we get \(\int_{\frac{3}{4}\pi}^{\infty}\frac{|\cos(t)|^{\sqrt{2}}}{t^3}\leq 0.0667\).
All in all, we have
\begin{align}
 \int_0^\infty\frac{\exp(-t^2/\sqrt{2})-|\cos(t)|^{\sqrt{2}}}{t^3}\mathrm{d}t
\leq 0.0312+0.2958-0.2577-0.0667=0.026\geq 0\notag
\end{align}
which is the last inequality to verify.
\subsection*{Acknowledgement} This article is an extension of the author's master thesis. The author likes
to thank his supervisor Hermann K\"onig for introducing him to this subject and for all advice and helpful discussions
during the preparation of this work.

\end{document}